\documentclass{elsart}
\journal{Theoretical Computer Science}

\newtheorem{lemme}{Lemma}

\usepackage{psfig}

\font\teuf=eufm10 scaled 1200
\font\seuf=eufm7 scaled 1200            
\font\sseuf=eufm5 scaled 1200

\newfam\euffam
        \textfont\euffam=\teuf
        \scriptfont\euffam=\seuf
        \scriptscriptfont\euffam=\sseuf

\newtheorem{theo}{Theorem}

\def\Choose#1#2{{{#1} \choose {#2}}}

\def\S0{ \{ S_0 \} }

\def\GC{{\mathcal G}}

\def\TC{{\mathcal T}}

\def\vect#1{\mbox{\boldmath{$#1$}}}

\font\tenbb=bbold10 scaled \magstep1

\font\sevenbb=bbold7 scaled \magstep1

\font\fivebb=bbold5 scaled \magstep1

\newfam\bbfam
\textfont\bbfam=\tenbb
\scriptfont\bbfam=\sevenbb
\scriptscriptfont\bbfam=\fivebb
\def\bb{\fam\bbfam\tenbb}

\def\1b{\mbox{$\bb 1$}}

\begin{document}

\begin{frontmatter}

\title{A formula for the number of tilings of an octagon by rhombi}

\author{N. Destainville$^{a}$,}
\author{R. Mosseri$^{b}$,}
\author{F. Bailly$^{c}$} 

\address{$^{a}$Laboratoire de Physique Th\'eorique~-- IRSAMC \\ 
UMR CNRS/UPS 5152,
Universit\'e Paul Sabatier, \\ 118, route de Narbonne, 31062
Toulouse Cedex 04, France.}

\address{$^{b}$Groupe de Physique des Solides, \\
CNRS et Universit\'e Paris 6, Campus Boucicaut, \\
140 Rue de Lourmel, 75015 Paris, France.}

\address{$^{c}$Laboratoire de Physique du Solide~-- CNRS, \\
1, place Aristide Briand, 92195 Meudon Cedex, France. }

\begin{abstract}

We propose the first algebraic determinantal formula to enumerate
tilings of a centro-symmetric octagon of any size by rhombi. This
result uses the Gessel-Viennot technique and generalizes to any
octagon a formula given by Elnitsky in a special case.
\end{abstract}

\end{frontmatter}

\section*{Introduction}

\begin{figure}[ht]
\begin{center}
\ \psfig{figure=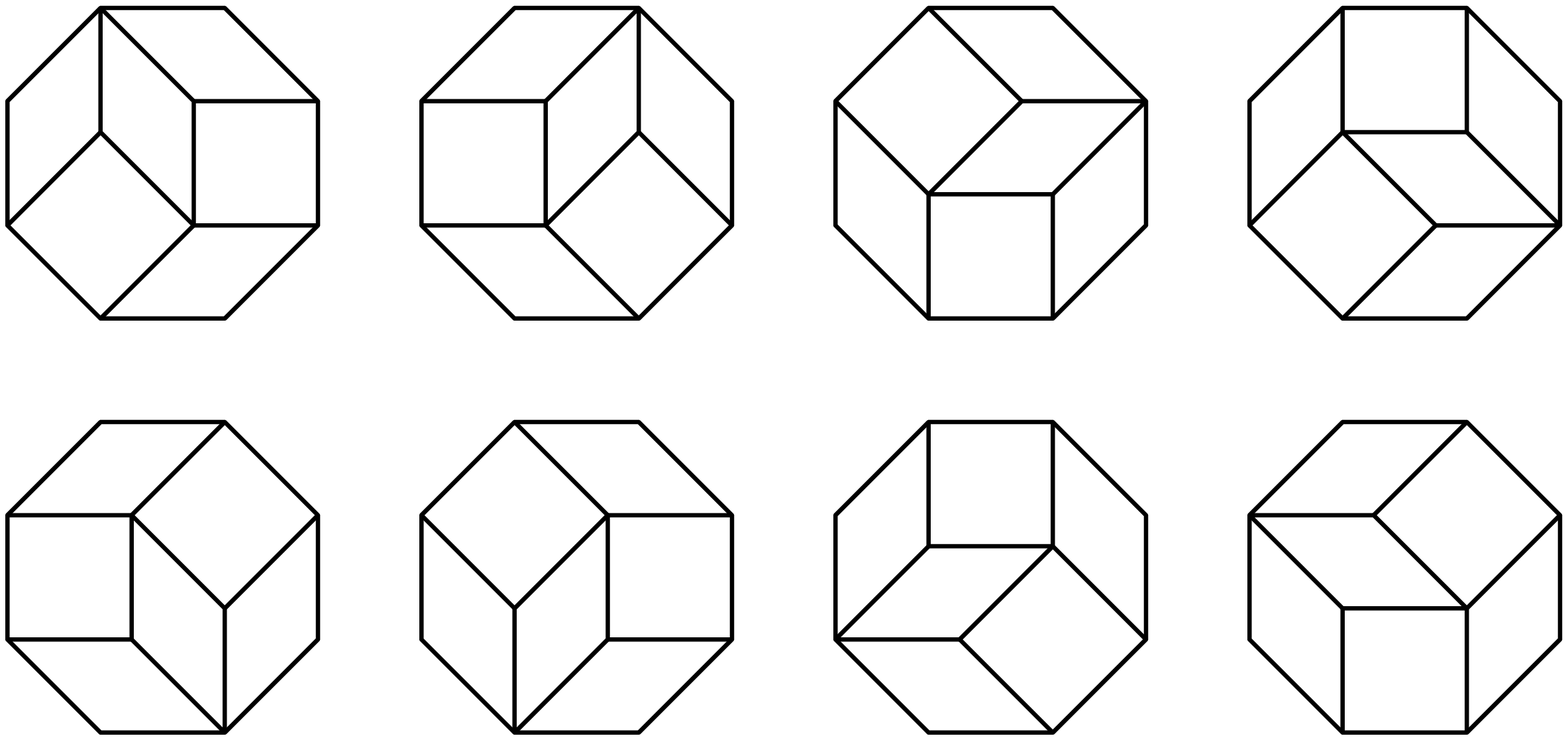,width=9cm} \
\end{center}
\caption{The eight tilings of the set $\TC_{1,1,1,1}$.}
\label{8octos}
\end{figure}

\begin{figure}[ht]
\begin{center}
\ \psfig{figure=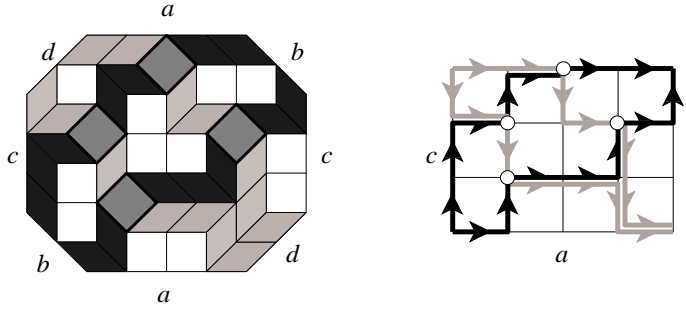,width=9cm} \
\end{center}
\caption{Left: Example of octagonal tiling of a centro-symmetric
octagon of sides $a,b,c,d$. There are 6 species of tiles with unitary
side length: two squares and four 45$^{\mbox{\scriptsize o }}$
rhombi. One species of squares, the ``tilted'' squares, is emphasized
in medium gray. They lie at the intersections of the light gray and
dark gray de Bruijn lines. De Bruijn lines are defined in
section~\ref{octogrid} and are made up of adjacent rhombi sharing an
edge with a given orientation.  Right: square grid representation of
the same tiling, obtained by shrinking all colored tiles; see
section~\ref{octogrid}. The white disks on the grid keep track of the
position of the tilted squares of the original tiling. The orientation
of edges (arrows) will be discussed later.}
\label{pavgrid}
\end{figure}

The enumeration of tilings of a centro-symmetric polygon by rhombi is
a notoriously difficult problem that concerns discrete mathematics and
theoretical computer science, as well as theoretical physics, in
relation with quasicrystallography.  In the latter community, these
tilings are usually called ``random tilings with octagonal symmetry''.
We address the following issue: given a centro-symmetric octagon
$O_{a,b,c,d}$, of integral sides lengths $a,b,c$ and $d$ (read
clockwise; see figure~\ref{pavgrid}, left), in how many ways is it
possible to fill it entirely, without any gap or overlap, with the
following six species of tiles: two differently oriented squares, and
four differently oriented 45$^{\mbox{\scriptsize o }}$ rhombi, the six
of them with unitary side lengths? So far, this question has been
solved in very particular instances only. We denote by $\TC_{a,b,c,d}$
the set of all the tilings of $O_{a,b,c,d}$, and by
$T_{a,b,c,d}=|\TC_{a,b,c,d}|$ the cardinality of $\TC_{a,b,c,d}$. For
example, figure~\ref{8octos} displays the eight tilings of the set
$\TC_{1,1,1,1}$.

Small systems have been studied in references~\cite{Mosseri93B,Bibi01}
up to sizes of some hundred tiles (see table~\ref{enum}). However, the
technique employed cannot reasonably provide tiling enumerations for
bigger octagons.
\begin{table}[ht]
\caption{\label{enum}Some tiling enumerations computed in
ref.~\cite{Bibi01}. The number of rhombi is given in column 3.}
\begin{center}
\begin{tabular}{|c|l|l|}
\hline
$a,b,c,d$ & $T_{a,b,c,d}$ & \mbox{\# tiles} \\ \hline \hline
1,1,1,1 & 8 \hfill (see fig.~\ref{8octos}) & 6 \\
2,2,2,2 & 5383 & 24 \\
3,3,3,3 & 273976272 & 54 \\
4,4,4,4 & 1043065776718923 & 96 \\
5,5,5,5 & 296755610108278480324496 & 150 \\
\hline
\end{tabular}
\end{center}
\end{table}
On the other hand, Elnitsky gave in ref.~\cite{Elnitsky97} two
formulas when two sides of the octagon are set to 1:
\begin{equation}
T_{a,1,c,1}=\sum_{r+s=a} \ \sum_{t+u=c} \Choose{r+t}{r}
\Choose{s+t}{s} \Choose{r+u}{r} \Choose{s+u}{s},
\label{Elnit1}
\end{equation}
and
\begin{equation}
T_{a,b,1,1}={2 (a+b+1)! \; (a+b+2)!
\over a! \; b! \; (a+2)! \; (b+2)!}.
\label{Elnit2}
\end{equation}
The first formula has been later partially simplified~\cite{Bibi01}:
\begin{equation}
T_{a,1,c,1} =  {(a+c+1)! \over a! \; c! \; (2a+1) (2c+1)} 
\left[ {2(a+c+1)! \over a! \; c!} 
+ \sum_{k=0}^a {1 \over 2k - 1} \Choose{a}{k} \Choose{c}{k} \right]
\label{francis}
\end{equation}
where the last sum can be written in terms of a hypergeometric
function 
\begin{equation}
\sum_{k=0}^a {1 \over 2k - 1} \Choose{a}{k} \Choose{c}{k} = \
_3\!\, F_2 \left[ -1/2,-a,-c ; 1/2,1 ;
1 \right].
\end{equation}

We propose a generalization of the first formula~(\ref{Elnit1})
to any side lengths,
where $T_{a,b,c,d}$ is written as a sum of products of determinants.
Even if the complexity of our formula increases with the system size,
it is the first explicit algebraic expression to count tilings
of an octagon (see eq.~(\ref{W})), which can be in principle calculated
for any system size. 

As it is discussed below into detail, Elnitsky's proof uses a ``square
grid representation'' of tilings (see figure~\ref{pavgrid}), which is
closely related to the ``de Bruijn dualization'', a wide-spread
technique in quasicrystal science. This dualization has proved
powerful to handle rhombus tilings in several
circumstances. The present paper also uses this technique, thus
generalizing Elnitsky's proof.

We now state our main result. Given the side lengths $a$, $b$, $c$ and $d$,
we denote by $X$ (resp. $Y$) the set 
of families of integers $(x_{k,l})$ (resp. $(y_{k,l})$), $k=1,\ldots,b$,
$l=1,\ldots,d$, satisfying the relations:
\begin{eqnarray}
\label{C1}
0 \leq x_{k,l} \leq a, & \phantom{\mbox{ and }} & \\
\label{c1}
x_{k,l} \leq x_{k',l'}& \mbox{ if } & k\leq k' \mbox{ and } l\leq l', \\
\label{C2}
0 \leq y_{k,l} \leq c, & \phantom{\mbox{ and }} & \\
\label{c2}
y_{k,l} \leq y_{k',l'}& \mbox{ if } & k\geq k' \mbox{ and } l\leq l'. 
\end{eqnarray} 
Note that conditions (\ref{c1}) and (\ref{c2}) are not exactly
similar. In the following, these integers will be the coordinates of
the white disks in the grid representation (see figure~\ref{pavgrid},
right, and section \ref{octogrid}). In addition, we set by convention
for $k=1,\ldots,b$ and $l=1,\ldots,d$
\begin{equation}
\label{conventions}
\begin{array}{ll}
x_{k,0}=0, & y_{k,0}=0, \\
x_{0,l}=0, & y_{0,l}=c, \\
x_{b+1,l}=a, & y_{b+1,l}=0, \\
x_{k,d+1}=a, & y_{k,d+1}=c.
\end{array}
\end{equation}
The reasons for this convention will be explicited below.
For any two such sequences $x=(x_{k,l})$ and $y=(y_{k,l})$,
we define the matrices $M^{(u)}(x,y)$ and $P^{(v)}(x,y)$ as follows:
$M^{(u)}(x,y)$ is a $b \times b$ matrix of coefficients
\begin{equation}
m_{ij} = \Choose{x_{j,u}-x_{i,u-1}+y_{j,u}-y_{i,u-1}}{x_{j,u}-x_{i,u-1}+j-i}
\end{equation}
for $1 \leq i,j \leq b$ 
and $P^{(v)}(x,y)$ is a $d \times d$ matrix of coefficients
\begin{equation}
p_{ij} = \Choose{x_{v,j}-x_{v-1,i}+y_{v-1,i}-y_{v,j}}{x_{v,j}-x_{v-1,i}+j-i}.
\end{equation}
for $1 \leq i,j \leq d$.
Note that, by convention, we set $\Choose{A}{B}=0$ whenever $A<0$ or
$B<0$ or $B>A$.

\bigskip

\begin{theo} 
\label{th1}
With the previous notations, the number of rhombus tilings
of a centro-symmetric octagon of sides $a,b,c,d$ reads:
\begin{equation}
\label{W}
T_{a,b,c,d} = \sum_{(x,y)\in X \times Y} \ \prod_{u=1}^{d+1} \det M^{(u)}(x,y)
\prod_{v=1}^{b+1} \det P^{(v)}(x,y).
\end{equation}
\end{theo}
It is demonstrated below that the determinants come from the
enumeration, by the Gessel-Viennot method (presented below), of
tilings of independent sub-domains of the octagon delimited by some points of
coordinates $(x_{k,l},y_{k,l})$ in the square grid representation. 

When $b=d=1$, the previous expression is reduced to Elsnitsky's
relation~(\ref{Elnit1}). Note that by contrast,
relation~(\ref{Elnit2}) is {\em not} a spacial case of this formula.
Beyond this simple case, the number of terms in the formula grows with
the octagon size. For example, for $(a,b,c,d)=(2,2,2,1)$, the formula
contains $6\times6=36$ terms to count the 480 tilings. For
$(a,b,c,d)=(2,2,2,2)$, there are $20 \times 20 =400$ terms and 5383
tilings. More generally, the number of terms grows exponentially with
the number of tiles, but it nevertheless grows exponentially more
slowly than the number of tilings. As a consequence, this formula is
exponentially more compact than the crude enumeration of tilings. This
point is discussed in the conclusion.

\section{Octagonal tilings and the square grid representation}
\label{octogrid}

\begin{figure}[ht]
\begin{center}
\ \psfig{figure=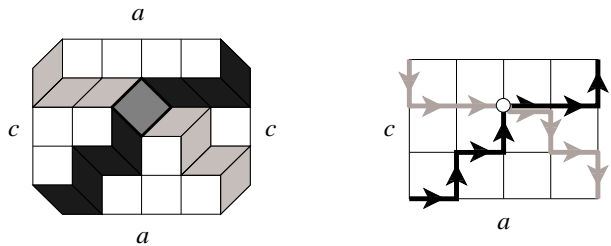,width=8cm} \
\end{center}
\caption{Example of tiling in the case $b=d=1$ (left) together with
its square grid representation (right). We have emphasized two de
Bruijn lines of the tiling, belonging to two different families. Their
intersection coincides with the gray tilted square.  They are
represented by two directed paths joining two opposite corners of the
grid. These paths can have multiple intersections.  To avoid ambiguity
on the position of the gray tilted square on the original tiling, a
white circle marks a {\em distinguished} vertex to keep track of its
position on the grid. Note that when two paths are tangent, we have
slightly shifted one of them for the sake of readability. However,
these paths are in fact superimposed and lie on the same grid edge.}
\label{bd1}
\end{figure}

In this section, we show that octagonal tilings are
conveniently represented by families of directed paths running on a
rectangular patch of square grid. This representation was used by
Elnitsky \cite{Elnitsky97} and it is reminiscent of the prior ``de
Bruijn dualization''~\cite{DeBruijn1,DeBruijn2,Socolar,Gahler} and derived
representations~\cite{Mosseri93B,Bibi01}.  We first expose briefly the
de Bruijn dualization process. Figures~\ref{pavgrid} and~\ref{bd1}
will help the reader. To begin with, we notice that tile edges can
have four possible orientations. We define a {\em family of de Bruijn
lines} for each orientation: de Bruijn lines are made up of adjacent
rhombi sharing an edge with this orientation. It is always possible to
extend these lines through the whole tiling up to a boundary tile. Two
examples of lines are presented in both figures~\ref{pavgrid} and
\ref{bd1}, belonging to two different families. A rhombic tile is
situated at the intersection of two lines of different families (see
the figures). By construction, lines of a same family never intersect;
there are respectively $a$, $b$, $c$ and $d$ lines in each family.

Now we show how to translate the de Bruijn's representation of a
tiling into its square grid representation. In figure~\ref{bd1}, we
show this correspondence in the simplest case
$b=d=1$~\cite{Elnitsky97}. The idea is to shrink the de Bruijn lines
of two families among four, so that they become paths on a square
grid, as displayed in figures~\ref{pavgrid} or \ref{bd1}. Because all
tiles of a de Bruijn line have an edge of a given orientation, these
paths are directed. The $b$ paths of the first family (denoted by
$SW$) go from the south-west corner to the north-east one (dark gray);
they can follow eastward and northward edges only; the $d$ paths of
the second family (denoted by $NW$) go from the north-west corner to
the south-east one (light gray); they can follow eastward and
southward edges only.

In the simplest case $b=d=1$, to avoid ambiguity due to path tangency
and to make this correspondence bijective~\cite{Elnitsky97}, we keep
track of the intersection of the de Bruijn lines thanks to a {\em
distinguished vertex}, represented by a white disk in the right
figure. It marks the position of the unique tilted square (medium
gray).  When $b>1$ or $d>1$ as in figure~\ref{pavgrid}, there are $bd$
intersections and therefore $bd$ tilted squares. Each of them must be
located by a distinguished vertex on the square grid. Paths do not
cross in a same family even though {\em they can be locally adjacent}
(see figure~\ref{pavgrid}).  We denote the paths of $SW$ (resp. $NW$)
by $SW_1,\ldots,SW_b$ (resp. $NW_1,\ldots,NW_d$) from left to right.
As a consequence, distinguished
vertices are indexed by two integers $k$ and $l$, and are denoted by
$DV_{k,l}$ (see figure~\ref{patches}).

\section{The Gessel-Viennot method}
\label{GV}

The Gessel-Viennot method~\cite{Gessel85,Stembridge90} is a
combinatorial technique for the counting of configurations of directed
non-intersecting paths on oriented graphs. This technique has already
proved very useful for the enumeration of rhombus tilings (see
\cite{Elnitsky97,Bibi98} for examples, as well as
section~\ref{proofth1} of the present paper). It has been extensively
described in the literature~\cite{Gessel85,Stembridge90} and we shall
only briefly explain it in the present paper, focusing on the
underlying ideas and not on technical details. The method is rather
general and can be applied to any acyclic oriented graph $\GC$, in
which are selected two families of vertices, $d_i$ (``departure''
vertices) and $a_j$ (``arrival'' vertices), $i,j=1,\ldots,n$. We
consider directed paths, running on $\GC$, starting from one vertex
$d_i$ and arriving at one vertex $a_j$. By ``directed'', we naturally
mean that the paths must follow the edge orientations. In addition,
this graph is supposed to satisfy the property of {\sl compatibility}:
if two directed paths on $\GC$ are going respectively from $d_{i_1}$
to $a_{j_1}$ and from $d_{i_2}$ to $a_{j_2}$, if these paths do not
cross, and if $i_1 < i_2$ then $j_1 < j_2$. This property is
very specific to two-dimensional graphs.

We are interested in the number $D_n$ of configurations of $n$
non-intersecting directed paths on $\GC$, where the $i$-th path goes
from $d_i$ to $a_i$: two paths are said to be non-intersecting if they
share no vertex; $n$ paths are said to be non-intersecting if any two
paths are non-intersecting. If we denote by $\lambda_{ij}$ the number
of paths going from $d_i$ to $a_j$, then the Gessel-Viennot theorem
states that
\begin{equation}
D_n=\det(\lambda_{ij})_{1 \leq i,j \leq n}.
\end{equation}

The idea of the proof is that in this determinant, all
configurations of $n$ paths, whether intersecting or not, the $i$-th
path going from $d_i$ to $a_{\sigma(i)}$, for any permutation
$\sigma$, are counted, with a $+$ or $-$ sign. Because of these
signs, all configurations with one or more intersections cancel two
by two. Only the non-intersecting configurations remain. They are
exactly the configurations under interest thanks to the property of
compatibility. The interested reader is referred to
Stembridge~\cite{Stembridge90} for more detailed explanations.

\section{Proof of theorem~\ref{th1}}
\label{proofth1}

We are now ready to prove theorem~\ref{th1}.  First of all, we need to
endow the square grid with integer coordinates in order to locate the
positions of distinguished vertices. They are defined according to the
usual conventions, so that the south-west and north-east corners have
respective coordinates $(0,0)$ and $(a,c)$. The coordinates of
$DV_{k,l}$ are denoted by $(x_{k,l},y_{k,l})$. By extension, we also
define the vertices $DV_{k,0}=(0,0)$, $DV_{k,d+1}=(a,c)$,
$DV_{0,l}=(0,c)$ and $DV_{b+1,l}=(a,0)$. They are the ends of paths of
families $SW$ and $NW$. If their coordinates are also denoted by
$(x_{k,l},y_{k,l})$, these last definitions are compatible with the
conventions~(\ref{conventions}) given in introduction.  All these
coordinates naturally obey relations~(\ref{C1}) and (\ref{C2}).

Furthermore, because of the directed character of de Bruijn lines,
distinguished vertices $DV_{k,l}$ are constrained by some conditions
when they belong to the same paths, and they must obey
relations~(\ref{c1}) and (\ref{c2}) as well. These four conditions
define the sets $X$ and $Y$, as it was stated in the introductory part.

Let $x=(x_{k,l}) \in X$ and $y=(y_{k,l}) \in Y$ be an admissible set
of coordinates of the distinguished vertices. We denote by $\TC_{x,y}$
the subset of tilings of $\TC_{a,b,c,d}$ in the square grid
representations of which the distinguished vertices have these
coordinates. The subsets $\TC_{x,y}$ are two-by-two disjoint so that
$T_{a,b,c,d} = \sum_{(x,y) \in X \times Y} |\TC_{x,y}|$.  Our purpose
is now to calculate each $|\TC_{x,y}|$.  This calculation is feasible
because for a given $(x,y)$, the subset $\TC_{x,y}$ can be factorized
into simple sets (see eq.~(\ref{factorization})). Each of them can in
turn be counted by the Gessel-Viennot method, which leads to
relation~(\ref{W}).

\begin{figure}[ht]
\begin{center}
\ \psfig{figure=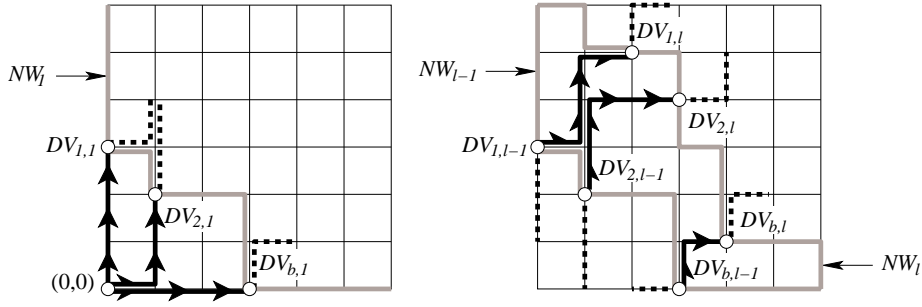,height=4cm} \
\end{center}
\caption{Examples of configurations of non-crossing directed paths,
belonging to the sets $sw(1)$ (left) and $sw(l)$ (right). Each path
(solid black lines) is a section of a path $SW_k$ and goes from
$DV_{k,l-1}$ to $DV_{k,l}$ (white circles).  By convention,
$DV_{k,0}=(0,0)$. The gray paths $NW_1$, $NW_{l-1}$ and $NW_l$ do not
belong to the path configurations and are drawn for information
only. The dashed black lines represent the possible continuations of
the original paths $SW_k$.  In this figure also, we have slightly
shifted tangent paths for the sake of readability, whereas in reality
they run on the same grid edge. }
\label{patches}
\end{figure}

First we need to introduce two definitions. Given a configuration of
vertices $DV_{k,l}=(x_{k,l},y_{k,l})$, we fix $l$ and we consider in
isolation the vertices $DV_{k,l-1}$ as well as $DV_{k,l}$,
$k=1,\ldots,b$ (see figures~\ref{patches} and
\ref{recomposition}). Then we define the set $sw(l)$ of {\em all} the
configurations of $b$ directed non-crossing paths, the $k$-th path
going from $DV_{k,l-1}$ to $DV_{k,l}$, with $k=1,\ldots,b$. These
paths are directed from south-west to north-east. They have no
constraint except that they are directed and non-crossing (these paths
can have tangencies).  In a similar way, we define the sets $nw(k)$
for any $k$: they are the sets of all configurations of $d$ directed
non-crossing paths, going from north-west to south-east.  The $l$-th
path goes from $DV_{k-1,l}$ to $DV_{k,l}$.

Now in order to prove theorem~\ref{th1}, we start from the following
observation, illustrated in figure~\ref{patches}: in $\TC_{x,y}$ the
distinguished vertices $DV_{k,l}$ are held fixed and one can consider
{\em independent} patches of the families $SW$ or $NW$, as
follows. Without loss of generality, we focus on $SW$. We cut each
path $SW_k$ into $d+1$ sections, denoted by $SW_k(l)$, where
$l=1,\ldots,d+1$. The section $SW_k(l)$ goes from $DV_{k,l-1}$ to
$DV_{k,l}$. Then all the $b$ sections $SW_k(l)$, $k=1,\ldots,b$ form a
local path configuration denoted by $P_{SW}(l)$. It belongs to
$sw(l)$.  In a similar way, the corresponding grid patches defined
with respect to the family $NW$ are denoted by $P_{NW}(k)$ and belong
to $nw(k)$.  Therefore, when $x$ and $y$ are fixed, we have the
natural inclusion
\begin{equation}
\TC_{x,y} \subset \prod_{u=1}^{d+1} sw(u) \prod_{v=1}^{b+1} nw(v),
\end{equation}
where the products are direct. We prove below that
\begin{lemme}
\label{lemme1}
The previous inclusion is an equality:
\begin{equation}
\label{factorization}
\TC_{x,y} = \prod_{u=1}^{d+1} sw(u) \prod_{v=1}^{b+1} nw(v).
\end{equation}
\end{lemme}

\medskip

\noindent It follows from eq.~(\ref{factorization}) that
\begin{equation}
\label{decomp}
\TC_{a,c,b,d} = \bigcup^{\mbox{\scriptsize disjoint}}_{(x,y) \in X \times Y}
\ \prod_{u=1}^{d+1} sw(u)
\prod_{v=1}^{b+1} nw(v),
\end{equation}
and that
\begin{equation}
T_{a,c,b,d} = \sum_{(x,y)\in X \times Y}
\ \prod_{u=1}^{d+1} |sw(u)|
\prod_{v=1}^{b+1} |nw(v)|.
\end{equation}
The remainder of the proof consists in calculating the cardinalities
$|sw(u)|$ and $|nw(v)|$ by the Gessel-Viennot method. Indeed, it is
also demonstrated below that
\begin{lemme}
\label{lemme2}
When $x$ and $y$ are fixed, 
\begin{equation}
|sw(u)| = \det M^{(u)}(x,y) \ \mbox{;}
\qquad |nw(v)| = \det P^{(v)}(x,y) .
\end{equation}
\end{lemme}

\subsection*{Proof of lemma~\ref{lemme1}:} 
\begin{figure}[ht]
\begin{center}
\ \psfig{figure=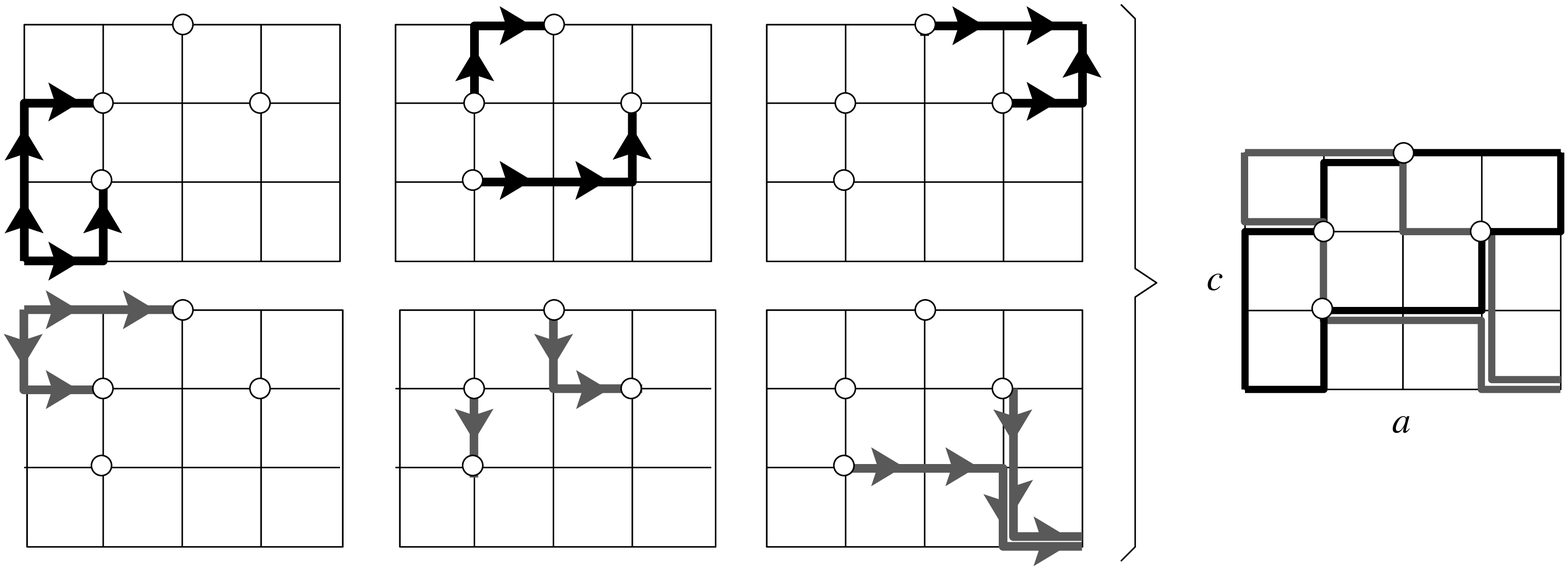,width=12cm} \
\end{center}
\caption{Example of path concatenation involved in the proof of
lemma~\ref{lemme1}.  The first line displays configurations from the
sets $sw(u)$ (black paths). The second line from the sets $nw(v)$
(gray paths). Altogether, the superposition of the six patches forms a
square grid representation.  Black paths and gray paths can cross only
at the distinguished vertices (white circles).}
\label{recomposition}
\end{figure}

We need to prove the
reverse inclusion
\begin{equation}
\TC_{x,y} \supset \prod_{u=1}^{d+1} sw(u) \prod_{v=1}^{b+1} nw(v).
\end{equation}
Configurations from the sets $sw(u)$ provide sections of paths from
$DV_{k,l-1}$ to $DV_{k,l}$. When concatenated, these sections provide
complete directed non-crossing paths from $(0,0)$ to $(a,c)$, which
form a family $SW$. In a similar way, sections from the $nw(v)$
provide directed non-crossing paths from $(0,c)$ to $(a,0)$, forming a
family $NW$. We only need to check that any two paths from $SW$ and
$NW$ only cross at the distinguished vertices $DV_{k,l}$. This point
is ensured by the directed character of path sections (see
figure~\ref{recomposition}). This last observation is crucial and all the
demonstration relies on it: it ensures the reverse inclusion and
therefore the direct character of the product~(\ref{factorization}),
from which our enumerating formula ensues.

\subsection*{Proof of lemma~\ref{lemme2}:}

So far we have used the terminologies ``non-{\em intersecting} paths''
in section~\ref{GV} and ``non-{\em crossing} paths'' in
section~\ref{proofth1}. Now it is time to precise what aspects these
two terms cover. We have seen that non-crossing paths can have
tangencies, that is to say they can share vertices or edges of
the grid, but they cannot step over one another. In particular,
non-crossing paths of families $SW$ (or $NW$) share their ends,
but can be indexed from west to east without ambiguity.

\begin{figure}[ht]
\begin{center}
\ \psfig{figure=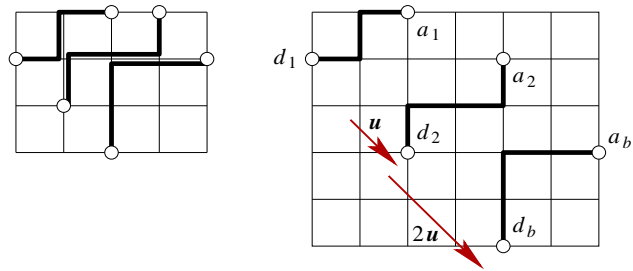,width=8.3cm} \
\end{center}
\caption{Left: a configuration of non-{\em crossing} paths; Right: 
a configuration of non-{\em intersecting} paths. The second 
configuration is obtained from the first one {\em via} the shifts
defined in the text. This correspondence is one-to-one.}
\label{shift}
\end{figure}

On the contrary, non-intersecting paths cannot share any vertex or
edge. Therefore, if we want to use the Gessel-Viennot method, we need
to transform configurations of non-crossing paths on the square grid
into configurations of non-intersecting paths.  The trick consists in
shifting non-crossing paths, as it is illustrated in
figure~\ref{shift}. The trick is standard and was already used
by Elnitsky \cite{Elnitsky97} for example.

We use a unitary shift vector $\vect{u}=(1,-1)$ and we shift the
$k$-th path section $SW_k(l)$ by a vector $(k-1) \vect{u}$ (see the
figure). The new paths still belong to the square grid. The $k$-th new
path section goes from the new vertex $DV'_{k,l-1}$ to the new vertex
$DV'_{k,l}$, of coordinates
\begin{equation}
x'_{k,l}=x_{k,l}+(k-1) \qquad \mbox{and} \qquad y'_{k,l}=y_{k,l}-(k-1).
\end{equation}
The so-obtained path configuration is non-intersecting by construction. 
This correspondence between configurations of non-crossing paths and 
configurations of non-intersecting paths is bijective.

Now we use the Gessel-Viennot technique by setting the departure
vertices $d_i=DV'_{i,l-1}$ and the arrival ones $a_j=DV'_{j,l}$.  The
number of directed paths on the square grid going from a vertex
$D=(x_D,y_D)$ to a vertex $A=(x_A,y_A)$ is simply given by the binomial
coefficient
\begin{equation}
\lambda = \Choose{x_A-x_D+y_A-y_D}{x_A-x_D}.
\end{equation}
Then one computes the coefficients $\lambda_{ij}$ involved in the
Gessel-Viennot method: $\lambda_{ij}=m_{ij}$. Thus one obtains the
matrices $M^{(u)}(x,y)$, the determinants of which count the
elements of $sw(u)$. In the same way, to count the elements of the
sets $nw(v)$, one must shift the sections of paths $NW_{l}$ by
$(l-1)\vect{v}$ where $\vect{v}=(1,1)$. One gets the matrices
$P^{(v)}(x,y)$ and takes their determinant, which completes the proof
of lemma~\ref{lemme2}.

\section{Conclusion}

We have demonstrated how Elnitsky's technique can be generalized 
to octagons of any size, leading to an explicit enumerative
formula (theorem~\ref{th1}).

We also notice that the conditions~(\ref{C1}) and (\ref{c1}) (resp.
(\ref{C2}) and (\ref{c2})) that define the set $X$ (resp. $Y$) are
identical to the conditions defining plane partitions of height $a$
(resp. $c$) on a $b \times d$ grid~\cite{MacMahon}. This point is
remarkable because such plane partitions are known to be equivalent to
rhombus tilings filling a centro-symmetric {\em hexagon} of
sides lengths $b$, $d$ and $c$ (resp. $b$, $d$ and $a$)~\cite{Elser}.
We have derived a partial combinatorial interpretation of our
formula~(\ref{W}) in terms of these tilings of hexagons. It is related
to a natural decomposition of the configuration sets of tilings of
octagons, as described in ref.~\cite{Bibi01}. But it goes beyond the
scope of the present paper and will be described
elsewhere~\cite{Prep}. 

If $T^{\mbox{\scriptsize hex}}_{a,b,c}$ denotes the number of tilings
of the centro-symmetric hexagon of sides $a$, $b$ and $c$, the
previous remark leads to the lower bound
\begin{equation}
T_{a,b,c,d} \geq T^{\mbox{\scriptsize hex}}_{b,d,c} \
T^{\mbox{\scriptsize hex}}_{b,d,a},
\end{equation}
the number of terms the formula. By construction, the sets
$nw(k)$ and $sw(l)$ are not empty and all terms are positive.
In statistical physics and more
specifically in quasicrystal science, people are interested in
thermodynamic quantities such as the configurational entropy (per
tile): $S=\ln(T_{a,b,c,d})/N_T$ where $N_T$ is the number of
tiles. With our polygonal boundary conditions, this quantity has a
finite limit when $N_T$ goes to infinity provided the relative ratios
of the side lengths also have a finite
limit~\cite{Mosseri93B,Bibi01,Bibi98,Elser}. In the so-called
``diagonal'' case where all side lengths are equal, taking into
account the number of tiles, the previous relation becomes $S \geq
S^{\mbox{\scriptsize hex}} = (3/2)\ln 3 - 2 \ln 2 \simeq
0.262$~\cite{Elser}.  The actual value of $S$ is numerically known to
be close to 0.36~\cite{Bibi01}. The previous lower bound is manifestly
loose and its improvement requires a better knowledge of the
asymptotic behavior of the determinants in~(\ref{W}) at the large size
limit.

But the main advantage of our formula precisely lies on the fact that
the previous bound is weak: the formula realizes an exponential
reduction of the number of terms as compared to a crude enumeration of
tilings. Indeed, as it was just discussed in the previous paragraph,
the number of terms grows exponentially like $\exp(0.26 N_T)$ wheres
the number of tilings grows like $\exp(0.36 N_T)$.  Even if in
practice we cannot compute numerically the number of tilings of
octagons bigger than in table~\ref{enum}, the progress is already
significant. Moreover there exists some hope to simplify our formula,
at least partially, as in eq.~(\ref{francis}).

In addition, the formula brings a new insight
into the structure of tiling sets: it emphasizes a natural
decomposition of the sets into smaller disjoint subsets, the cardinality of
which is simply given by evaluation of determinants.

\section*{Acknowledgments} 

We are indebted to one of the referees for a careful reading of our
manuscript and sound suggestions of improvements.

\end{document}